\documentclass{amsart}

\newtheorem{theorem}{Theorem}[section]
\newtheorem*{theorem A}{Theorem A}
\newtheorem*{theorem B}{N\"olker's Theorem}

\newtheorem{result}{Result}[section]

\theoremstyle{remark}

\theoremstyle{remark}

\theoremstyle{definition}

\newtheorem{definition}{Definition}[section]

\numberwithin{equation}{section}
\def\({\left ( }
\def\){\right )}
\def\<{\left < }
\def\>{\right >}

\begin{document}

\vspace{2cm}

\title{Helices in the Euclidean 5-Space $ E^5$ }

\author{MELEK MASAL, A. ZEYNEP AZAK}
\address{ Sakarya University, Faculty of Education, Department of Elementary Education\\
    Hendek- Sakarya-TURKEY$^1$\\}

\email{mmasal@sakarya.edu.tr,apirdal@sakarya.edu.tr}



\begin{abstract}
In this study, we have identified $V_3$ slant helix ($2^{nd}$ type
slant helix, $V_5$ slant helix ($3^{rd}$ type slant helix)  and
attained some characteristic properties in the Euclidean 5-Space
$E^5$. In addition to this, we have proven that there are no other
helices other than $V_1$ helix (inclined curve), $V_3$ slant helix
and $V_5$ slant helix in the 5-dimensional Euclidean space $E^5$.
\end{abstract}
\maketitle
\section{Introduction}
\rm The theory of curves is one of the fundamental topics of
differential geometry. Some specific curves play important roles
in a variation of sciences. As an example, the helix curves can
often be seen in the fields of biology and computer technologies
along with the daily life, \cite{Al2}. The classic results in
$R^3$ that are related to the helix curves and have so many fields
of application, were given by M.A. Lancret in 1802 and by B. de
Saint Venant in 1845, \cite{Br}. There have been many studies
related to the slant helices and darboux helices in the Euclidean
3-Space,  \cite{Al2,Iz,Ku,Zp} and some results have been achieved
that are related to the helices and $B_2$ slant helix (3rd type
slant helix)  in the Euclidean 4-Space, \cite{Ka,Tr1,Tr2}. Apart
from that, while different characterizations have been given for
the inclined and non-null inclined curves in the Euclidean 5-Space
and Lorentzian space, \cite{Al1,Ba}. Moreover the non-null helices
have been examined in the Lorentzian 6-Space and new
characterizations have been reached for the  $V_n$ slant helix in
the $n$-dimensional Euclidean Space, \cite{Go,Bo}.\\In this study,
$V_3$ slant helix and $V_5$ slant helix have been identified and
some results have been obtained in the five dimensional Euclidean
space $E^5$. Then we have proven that there are no other helices
other than $V_1$-helix, $V_3$ slant helix and $V_5$ slant helix in
$E^5$.

\section{Preliminaries}Let  $\alpha {\rm :I} \subset {\rm R} \to {\rm E}^{\rm 5}$ be an arbitrary curve in
$E^5$. Recall that the curve $\alpha$ is said to be of unit speed
curve if $\langle \alpha '(s),\alpha '(s)\rangle  = 1$ where
$\langle ,\rangle$ is the standard scalar  product in the
Euclidean space $E^5$ given by
$$\langle X,Y\rangle  =\sum\limits_{i = 1}^5 {x_i } y_i$$\\for each $X = \left( {x_1 ,x_2 ,x_3 ,x_4 ,x_5 } \right),Y = \left(
{y_1 ,y_2 ,y_3 ,y_4 ,y_5 } \right) \in E^5$. In particular, the
norm of a vector $X$ is given by $\left\| X \right\| = \sqrt
{\left\langle {X,X} \right\rangle }$, \cite{Ha}.\\Let $\left\{
{V_1 ,V_2 ,V_3 ,V_4 ,V_5 } \right\}$ be the moving frame along
$\alpha$. The Frenet equations of the curve $\alpha$ are given by

$$\left[ \begin{array}{l}
 V'_1  \\
 V'_2  \\
 V'_3  \\
 V'_4  \\
 V'_5  \\
 \end{array} \right] = \left[ \begin{array}{l}
 \,\,\,0\,\,\,\,\,\,\,\,\,k_1 \,\,\,\,\,\,\,0\,\,\,\,\,\,\,\,0\,\,\,\,\,\,\,0 \\
  - k_1 \,\,\,\,\,\,\,0\,\,\,\,\,\,\,\,k_2 \,\,\,\,\,\,0\,\,\,\,\,\,\,0 \\
 \,\,\,0\,\,\,\,\, - k_2 \,\,\,\,\,\,0\,\,\,\,\,\,\,\,k_3 \,\,\,\,\,0 \\
 \,\,\,0\,\,\,\,\,\,\,\,\,0\,\,\,\, - k_3 \,\,\,\,\,\,0\,\,\,\,\,\,k_4  \\
 \,\,\,0\,\,\,\,\,\,\,\,\,0\,\,\,\,\,\,\,\,0\,\,\,\,\, - k_4 \,\,\,\,0 \\
 \end{array} \right]\left[ \begin{array}{l}
 V_1  \\
 V_2  \\
 V_3  \\
 V_4  \\
 V_5  \\
 \end{array} \right]$$
where $V_i$ are called $i^{th}$ Frenet vectors and the functions
$k_i$ are called $i^{th}$ curvatures of the curve $\alpha$,
\cite{Ha}.\\A regular curve is called a W-curve if it has constant
Frenet curvatures, \cite{Yl}.

\section{$V_1$ helices in the five dimensional Euclidean space $E^5$}

\begin{definition}
Let  $\alpha {\rm :I} \subset {\rm R} \to {\rm E}^{\rm 5}$ be a
unit speed curve. If the tangent vector $V_1$ of the curve
$\alpha$ makes a constant angle with the fixed direction $U$, then
$\alpha$ is either a $V_1$-helix or inclined curve, \cite{Al1}.
\end{definition}
Here $\left\langle {V_1 ,U} \right\rangle  = \cos \theta  =
cons\tan t \ne 0$ expression can be written from the definition
3.1.
\begin{theorem}
Let $\alpha$ be a unit speed regular curve in $E^5$. Then $\alpha$
is a $V_1$-helix if and only if the function,
$$\left( {\frac{{k_1 }}{{k_2 }}} \right)^2  + \frac{1}{{k_3 ^2
}}\left[ {\left( {\frac{{k_1 }}{{k_2 }}} \right)^\prime  }
\right]^2  + \frac{1}{{k_4 ^2 }}\left[ {\frac{{k_1 k_3 }}{{k_2 }}
+ \left[ {\frac{1}{{k_3 }}\left( {\frac{{k_1 }}{{k_2 }}}
\right)^\prime  } \right]^\prime  } \right]^2$$ is
constant.\\Furthermore;
$$U = \cos \theta \left[ {V_1  + \frac{{k_1 }}{{k_2 }}V_3  +
\frac{1}{{k_3 }}\left( {\frac{{k_1 }}{{k_2 }}} \right)^\prime  V_4
+ \frac{1}{{k_4 }}\left( {\frac{{k_1 k_3 }}{{k_2 }} + \left[
{\frac{1}{{k_3 }}\left( {\frac{{k_1 }}{{k_2 }}} \right)^\prime  }
\right]^\prime  } \right)V_5 } \right]$$ where $\theta$ is the
angle between the vectors $V_1$ and $U$, \cite{Al1}.
\end{theorem}
\begin{theorem}
Let $\alpha$ be a unit speed regular curve in $E^5$. Then,
$\alpha$ is a $V_1$-helix if and only if the following equalities
is satisfied
$$k_4 \,f\left( s \right) = \frac{{k_1 \,k_3 }}{{k_2 }} + \left[
{\frac{1}{{k_3 }}\left( {\frac{{k_1 }}{{k_2 }}} \right)^\prime  }
\right]^\prime$$
and
$$\frac{1}{{k_4 }}\frac{d}{{ds}}f\left( s \right) =  - \frac{1}{{k_3
}}\left( {\frac{{k_1 }}{{k_2 }}} \right)^\prime$$ where $f$ is
$C^2$-function, \cite{Al1}.
\end{theorem}
\begin{theorem}
$\alpha$ is a unit speed curve in $E^5$. Then $\alpha$ is a
$V_1$-helix if and only if the equation is satisfied.
$$\frac{1}{{k_3 }}\left( {\frac{{k_1 }}{{k_2 }}} \right)^\prime   =
\left( {A - \int {\left[ {\frac{{k_1 k_3 }}{{k_2 }}\sin \int {k_4
ds} } \right]ds} } \right)\sin \int {k_4 ds - \left( {B + \left[
{\int {\frac{{k_1 k_3 }}{{k_2 }}\cos \int {k_4 ds} } } \right]ds}
\right)} \cos \int {k_4 ds}$$  for some constant $A$ and $B$,
\cite{Al1}.
\end{theorem}
Now, we will examine the helices other than $V_1$-helix in the
five dimensional Euclidean space $E^5$.

\section{$V_3$ slant helix ($2^{nd}$ type slant helix) in the five dimensional Euclidean space  $E^5$}
\begin{definition}
Let $\alpha$ be a unit speed curve with nonzero curvatures
$k_i,(\,1 \le i \le 4)$ in $E^5$. If the third unit Frenet vector
field $V_3$ of the curve $\alpha$ makes a constant angle $\phi$
with the fixed direction $U$, then $\alpha$ is called  a
$V_3$-slant helix (or $2^{nd}$ type slant helix). (Suppose that
$\left\langle {U,U} \right\rangle  = 1$).
\end{definition}
\begin{theorem}
Let $\alpha$ be a unit speed  curve in $E^5$. Then $\alpha$ is a
$V_3$ slant helix if and only if $\frac{{k_2 }}{{k_1 }} = cons\tan
t$ and $\frac{{k_3 }}{{k_4 }} = cons\tan t$.
\end{theorem}
\textbf{Proof:} If $\alpha$ is a $V_3$-slant helix and $U$ is a
fixed unit vector, then the following equality can be written as
\begin{equation}
\langle V_3 ,U\rangle  = \cos \phi  = cons\tan t \ne 0.
\end{equation}
Taking the differential of equation (4.1) with respect to $s$ and
using the Frenet equations, we obtain
$$\langle  - k_2 V_2  + k_3 V_4 ,U\rangle  = 0.$$
Therefore, $U$ lies on the hyperplane spanned by the Frenet
vectors $V_1$,$V_3$ and $V_5$. Then, we reach
\begin{equation}
U = u_1 V_1  + u_3 V_3  + u_5 V_5
\end{equation}
where, $u_i  = u_i (s)$ and $u_3  = \cos \phi = cons\tan t$.\\The
differential of equation (4.2), we have
$$u'_1 \,V_1  + \left( {u_1 k_1  - u_3 k_2 }
\right)V_2  + u'_3 \,V_3  + \left( {u_3 k_3  - u_5 k_4 }
\right)V_4  + u'_5 \,V_5=0.$$ By the above equality, the
coefficients $V_i$ are zero for $1 \le i \le 5$. So we get
$$\begin{array}{l}
 u'_1  = 0 \\
 u_1 k_1  - u_3 k_2  = 0 \\
 u'_3  = 0 \\
 u_3 k_3  - u_5 k_4  = 0 \\
 u'_5  = 0 \\
 \end{array}$$
Thus, it is easy to obtain that the coefficients $u_1$,$u_3$ and
$u_5$ are given by
\begin{equation}
\begin{array}{l}
 u_1  = \cos \phi \frac{{k_2 }}{{k_1 }} = cons\tan t \\
 u_3  = \cos \phi  = cons\tan t \\
 u_5  = \cos \phi \frac{{k_3 }}{{k_4 }} = cons\tan t \\
 \end{array}
\end{equation}
Since the coefficients $u_1$ and $u_5$ constants, the ratios
$\frac{{k_2 }}{{k_1 }}$ and $\frac{{k_3 }}{{k_4 }}$ are constant,
respectively. Therefore, if we substitute $u_1$,$u_3$ and $u_5$ in
(4.2), we have
\begin{equation}
{\rm U = cos}\phi \,\frac{{k_2 }}{{k_1 }}\,V_1  + \cos \phi
\,\,V_3  + \cos \phi \,\frac{{k_3 }}{{k_4 }}\,V_5.
\end{equation}
where $\phi  \ne k\frac{\pi }{2}$.\\Conversely, while the ratios
$\frac{{k_2 }}{{k_1 }}$ and $\frac{{k_3 }}{{k_4 }}$ are constant,
we can define the vector $U$. Since the differential of $U$ is
$\frac{{dU}}{{ds}} = 0$, $U$ is a fixed vector. Furthermore, since
$\langle V_3 ,U\rangle  = \cos \phi  = cons\tan t$, the curve
$\alpha$ become $V_3$ slant helix.\\This completes the proof.
\begin{result}
$\alpha$ is a $V_3$ slant helix if and only if the ratios
$\frac{{\left\| {V'_2 } \right\|}}{{\left\| {V'_1 } \right\|}}$
and $\frac{{\left\| {V'_4 } \right\|}}{{\left\| {V'_5 }
\right\|}}$ are constant.
\end{result}
\textbf{Proof:} If $\alpha$ is a $V_3$ slant helix, then from the
theorem 4.1 we can write
$$\frac{{k_2 }}{{k_1 }} = cons\tan t$$
and
$$\frac{{k_3 }}{{k_4 }} = cons\tan t.$$
From the Frenet equations, it is easily to see that
$$\frac{{\left\| {V'_2 } \right\|}}{{\left\| {V'_1 } \right\|}} =
\sqrt {1 + \left( {\frac{{k_2 }}{{k_1 }}} \right)^2 }$$
and
$$\frac{{\left\| {V'_4 } \right\|}}{{\left\| {V'_5 } \right\|}} =
\sqrt {1 + \left( {\frac{{k_3 }}{{k_4 }}} \right)^2 }$$.\\So we
have the ratios
$$\frac{{\left\| {V'_2 } \right\|}}{{\left\| {V'_1 } \right\|}}$$ and $$\frac{{\left\| {V'_4 } \right\|}}{{\left\|
{V'_5 } \right\|}}$$ are constant.\\Conversely, if $\frac{{\left\|
{V'_2 } \right\|}}{{\left\| {V'_1 } \right\|}} = cons\tan t$ and
$\frac{{\left\| {V'_4 } \right\|}}{{\left\| {V'_5 } \right\|}} =
cons\tan t$, then the ratios $\frac{{k_2 }}{{k_1 }}$ and
$\frac{{k_3 }}{{k_4 }}$ are constant. Thus, from Teorem 4.1
$\alpha$ is $V_3$ slant helix.
\begin{result}
If $\alpha$ is a W-curve, then $\alpha$ is $V_3$ slant helix.
\end{result}
\textbf{Proof:} If $\alpha$ is a W-curve, then the curvatures $k_i
,\quad 1 \le i \le 4$ are constants. Thus, the ratios $\frac{{k_2
}}{{k_1 }}$ and $\frac{{k_3 }}{{k_4 }}$ constant. This shows that
$\alpha$ is $V_3$ slant helix from Theorem 4.1.

\section{$V_5$ slant helix ($3^{rd}$ type slant helix) in five dimensional Euclidean space $E^5$}
\begin{definition}
A unit speed curve $\alpha :I\subset R \to E^5$  is said to be a
$V_5$ slant helix ($3^{rd}$ type slant helix ) if the fifth unit
Frenet vector field $V_5$ makes a constant angle $\Psi$ with the
unit and fixed direction $U$.
\end{definition}
\begin{theorem}
Let $\alpha$  be a unit speed curve in $E^5$.\\ \textbf{i)}
$\alpha$ is a  $V_5$ slant helix if and only if the function
$$\frac{1}{{k_1 ^2 }}\left[ {\frac{{k_4 k_2 }}{{k_3 }} + f'\left( s
\right)} \right]^2  + \left[ {f\left( s \right)} \right]^2  +
\left( {\frac{{k_4 }}{{k_3 }}} \right)^2$$ is constant, where
$f(s) = \left( {\frac{{k_4 }}{{k_3 }}} \right)^\prime
\frac{1}{{k_2 }}$.\\ \textbf{ii)}$\alpha$ is a  $V_5$ slant helix
if and only if the following equation is satisfied
$$\left[ {\frac{{k_4 k_2 }}{{k_1 k_3 }} + \frac{{f'}}{{k_1 }}}
\right]^\prime   + f\,k_1  = 0$$ where $f(s) = \left( {\frac{{k_4
}}{{k_3 }}} \right)^\prime \frac{1}{{k_2 }}$.
\end{theorem}
\textbf{Proof:} \textbf{i)} From the above definition 5.1, we give
the following
\begin{equation}
\left\langle {V_5 ,U} \right\rangle  = \cos \psi  = cons\tan t
\end{equation}
If we take the differential of equation (5.1) with respect to $s$,
we obtain
$$\langle  - k_4 V_4 ,U\rangle  = 0.$$
Therefore, we may express
\begin{equation}
U = u_1 V_1  + u_2 V_2  + u_3 V_3  + u_5 V_5.
\end{equation}
we know that $u_i  = u_i \left( s \right)$ and $u_5  = \cos \psi =
cons\tan t$\\The differentiation (5.2) gives
$$\left( {u'_1  - u_2 k_1 } \right)V_1  + \left( {u_1 k_1  + u'_2  -
u_3 k_2 } \right)V_2  + \left( {u_2 k_2  + u'_3 } \right)V_3  +
\left( {u_3 k_3  - u_5 k_4 } \right)V_4  + u'_5 \,V_5  = 0$$ and
from this equation we find
$$\begin{array}{l}
 u'_1  - u_2 k_1  = 0 \\
 u_1 k_1  + u'_2  - u_3 k_2  = 0 \\
 u_2 k_2  + u'_3  = 0 \\
 u_3 k_3  - u_5 k_4  = 0 \\
 u'_5  = 0. \\
 \end{array}$$
Using the above equations, we can form
\begin{equation}
\begin{array}{l}
 u_1  = \frac{{\cos \psi }}{{k_1 }}\left\{ {\frac{{k_4 k_2 }}{{k_3 }} + \left[ {\left( {\frac{{k_4 }}{{k_3 }}} \right)^\prime  \frac{1}{{k_2 }}} \right]^\prime  } \right\} \\
 u_2  =  - \cos \psi \left( {\frac{{k_4 }}{{k_3 }}} \right)^\prime  \frac{1}{{k_2 }} \\
 u_3  = \cos \psi \frac{{k_4 }}{{k_3 }} \\
 u_5  = \cos \psi  = cons\tan t. \\
 \end{array}
\end{equation}
and
\begin{equation}
u'_1  = u_2 \,k_1.
\end{equation}
If we define $f = f(s)$ by
$$\left( {\frac{{k_4 }}{{k_3 }}} \right)^\prime  \frac{1}{{k_2 }} =
f\left( s \right)\,\,$$
then the equation (5.3) writes as
\begin{equation}
\begin{array}{l}
 u_1  = \frac{{\cos \psi }}{{k_1 }}\left( {\frac{{k_4 k_2 }}{{k_3 }} + f'} \right) \\
 u_2  =  - \cos \psi f \\
 u_3  = \cos \psi \frac{{k_4 }}{{k_3 }} \\
 u_5  = \cos \psi  = cons\tan t. \\
 \end{array}
\end{equation}
and equation (5.4) can be written as
\begin{equation}
\left[ {\frac{{k_4 k_2 }}{{k_1 k_3 }} + \frac{{f'}}{{k_1 }}}
\right]^\prime   + f\,k_1  = 0.
\end{equation}
Therefore, equation (5.2) takes the following form:
\begin{equation}
U = \cos \psi \left( {\frac{{k_4 k_2 }}{{k_1 k_3 }} +
\frac{{f'}}{{k_1 }}} \right)V_1  - \cos \psi f\,V_2  + \cos \psi
\frac{{k_4 }}{{k_3 }}V_3  + \cos \psi \,V_5.
\end{equation}
Since $U$ is a fixed vector, the following expression
\begin{equation}
\frac{1}{{k_1^2 }}\left[ {\frac{{k_2 k_4 }}{{k_3 }} + f'}
\right]^2  + f^2  + \left( {\frac{{k_4 }}{{k_3 }}} \right)^2
\end{equation}
is obtained as constant.\\Conversely, if equation (5.8) holds,
then the fixed vector $U$ can be defined as
$$U = \cos \psi \left( {\frac{{k_4 k_2 }}{{k_1 k_3 }} +
\frac{{f'}}{{k_1 }}} \right)V_1  - \cos \psi f\,V_2  + \cos \psi
\frac{{k_4 }}{{k_3 }}V_3  + \cos \psi \,V_5$$ In this case
$\frac{{dU}}{{ds}} = 0$ and $\langle V_5 ,U\rangle = \cos \psi  =
cons\tan
t$. It is clear that $\alpha$ is $V_5$ slant helix.\\
\textbf{ii)} If $\alpha$ is $V_5$ slant helix, then from the proof
of theorem 5.1 i), we have $$\left[ {\frac{{k_4 k_2 }}{{k_1 k_3 }}
+ \frac{{f'}}{{k_1 }}} \right]^\prime   + f\,k_1  = 0.$$
Conversely, if the equation (5.6) holds, then the following can be
written
$$U = \cos \psi \left( {\frac{{k_4 k_2 }}{{k_1 k_3 }} +
\frac{f}{{k_1 }}} \right)V_1  - \cos \psi f\,V_2  + \cos \psi
\frac{{k_4 }}{{k_3 }}V_3  + \cos \psi \,V_5.$$ Since
$\frac{{dU}}{{ds}} = 0$ and $\langle V_5 ,U\rangle  = \cos \psi  =
cons\tan t$, $\alpha$ becomes $V_5$ slant helix.
\begin{result}
 Let $\alpha$ be $V_5$ slant helix. If the ratio $\frac{{k_4 }}{{k_3
 }}$ is constant, then the ratio $\frac{{k_2 }}{{k_1 }}$ is
 constant. (Namely, $\alpha$ is $V_3$ slant helix)
\end{result}
\textbf{Proof:} Suppose that $\alpha$ is $V_5$ slant helix and the
ratio $\frac{{k_4 }}{{k_3 }}$ is constant. So the equation (5.8)
becomes constant, that is, the ratio  $\frac{{k_2 }}{{k_1 }}$ is
found as constant. This means that the curve $\alpha$ is $V_3$
slant helix from theorem 4.1.
\begin{result}
If $\alpha$ is $V_5$ slant helix and $\frac{{\left\| {V'_5 }
\right\|}}{{\left\| {V'_4 } \right\|}}$ is constant, then
$\frac{{\left\| {V'_2 } \right\|}}{{\left\| {V'_1 } \right\|}}$ is
constant. (Namely, $\alpha$ is $V_3$ slant helix)
\end{result}
\textbf{Proof:} It is obvious from result 5.1.

Also it can be examined if there are any other helices other than
$V_1$ helix, $V_3$ slant helix and $V_5$ slant helix in $E^5$.
\begin{theorem}
Let $\alpha$  be a unit speed curve in $E^5$.\\ \textbf{i)}There
is no fixed direction making a constant angle with the second
Frenet vector $V_2$ of the curve $\alpha$.\\ \textbf{ii)}There is
no fixed direction making a constant angle with the fourth Frenet
vector  $V_4$ of the curve $\alpha$.
\end{theorem}
\textbf{Proof:} Let us assume that the second Frenet vector $V_2$
of the unit speed curve $\alpha$ in $E^5$ makes a constant angle
with the fixed direction $U$. So, we can write
\begin{equation}
\langle V_2 ,U\rangle  = \cos \beta
\end{equation}
where $\beta$ is a constant angle between $V_2$ and $U$.
Differentiating equation (5.9) with respect to $s$, we obtain
$$\langle  - k_1 V_1  + k_2 V_3 ,U\rangle  = 0.$$
This shows that the vector $U$ is perpendicular to the Frenet
vectors $V_1$ and  $V_3$, so the following can be written
\begin{equation}
U = u_2 V_2  + u_4 V_4  + u_5 V_5 ,\quad \quad u_i  = u_i \left( s
\right)
\end{equation}
Differentiating the equation (5.10), we have
$$\left( { - u_2 k_1 } \right)V_1  +  u'_2 \,V_2  +
\left( {u_2 k_2  - u_4 k_3 } \right)V_3  + \left( {u'_4  - u_5 k_4
} \right)V_4  + \left( {u_4 k_4  + u'_5 } \right)V_5  = 0$$ which
leads to the following system
\begin{equation}
\begin{array}{l}
  - u_2 k_1  = 0 \\
 u'_2  = 0 \\
 u_2 k_2  - u_4 k_3  = 0 \\
 u'_4  - u_5 k_4  = 0 \\
 u_4 k_4  + u'_5  = 0 \\
 \end{array}
\end{equation}
By taking account of the equations (5.11) we get $u_2  = u_4  =
u_5  = 0$ which gives us that $\mathop U\limits^ \to   = \mathop
0\limits^ \to$.\\Moreover, this shows that there is no fixed
direction $U$ that makes a constant angle with the Frenet vector
$V_2$.\\ \textbf{ii)} The proof is similar to the proof of theorem
5.2.i.

\end{document}